\documentstyle[12pt,amstex,amssymb]{article}
\newcommand{\<}{\langle}
\renewcommand{\>}{\rangle}

\newcommand{\e}{\varepsilon}

\begin{document}
\setlength{\baselineskip}{18pt}

\begin{center}
{\Large\bf A new formulation of the
near-equilibrium}\smallskip\\
{\Large\bf theory of turbulence}
\footnote{Supported in part by the
Applied Mathematical Sciences subprogram of the Office of Energy
Research, U.S. Department of Energy, under contract
DE--AC03--76--SF00098, and in part by the National Science
Foundation under grants DMS94--14631 and DMS973-2710}

\vspace{.5 truein}

{\sc G.I.~Barenblatt and Alexandre J.~Chorin}

\bigskip
Department of Mathematics\\
University of California\\
Berkeley, California 94720--3840
\end{center}

\bigskip\bigskip\noindent
{\bf Abstract.}

We present  a status report on a discrete approach to the the
near-equilibrium
statistical
theory of three-dimensional turbulence,
which generalizes earlier work by no longer requiring
that the vorticity field be a union of discrete vortex filaments.
The idea is to take a special limit of a dense lattice vortex system,
in a way that brings out a connection between turbulence and critical
phenomena.
The approach produces statistics with basic features of turbulence, in
particular intermittency and
coherent  structures.
The numerical calculations have not yet been brought to convergence,
and  at present the results are only qualitative.

\newpage
\section{Introduction}

In recent work on wall-bounded turbulence
and local structure in fully developed turbulence [1,2,3,6] we have made
extensive use
of vanishing-viscosity asymptotics  and of expansions in powers of the
small parameter
$\frac{1}{\ln Re}$ ($Re$ is the Reynolds number) around  the asymptotic,
highly intermittent state $Re
\,\infty$; we concluded in particular
that the Kolmogorov-Obukhov spectrum must be
corrected when the viscosity is finite, and  also that moments of the
velocity field
of order $p\ge4$ diverge as the viscosity $\nu$ tends to zero.
Our analysis used results of a near-equilibrium statistical theory of
turbulence based on vortex representations [5,6].
In the present paper we
present briefly a new
derivation of the near-equilibrium  statistical theory of turbulence,
without assuming that the vorticity field can represented as a union of
discrete vortex filaments, and
discuss conclusions that can be drawn from it, in particular regarding
intermittency, scaling, and the existence of coherent structures.
The work includes preliminary numerical computations in three space
dimensions, and
while these computations are
informative, they are still
incomplete.

One of the keys to the understanding of turbulence is an explanation
of what can be viewed as a statistical equilibrium in fluid
mechanics. Statistical theories in both mechanics and physics are built along
differing lines when
they describe systems at equilibrium, near equilibrium, or far from
equilibrium.
The validity of the expansions we used in [1,2,3,6]
is connected with specific properties of statistical equilibria in fluid
mechanical vortex systems.
A fluid mechanical system
can be described macroscopically or microscopically;
a macroscopical description requires in principle just a few global
parameters, such as the
energy per unit mass or the rate of energy dissipation per unit mass;
A microscopical description purports to provide much of the detail in the
system. As many microscopical descriptions are generally compatible with a
single
macroscopical picture,  a microscopical description usually consists of an
ensemble of possible realizations,
all compatible with a fixed macroscopical picture, and each occurring with
a certain probability.
A statistical equilibrium is an ensemble of microscopical pictures,
each with
the probabilities it would have when the
macroscopical, large scale,  description has been kept constant for
a long time.
In other words, to find a statistical equilibrium, one fixes the large
scale parameters, and, having waited a while for the system to settle down,
one observes the resulting microscopic states
and the probabilities of their occurrence.
There are a number of recipes for constructing statistical equilibria.
One is the so-called ``microcanonical" ensemble, in which one
considers all the states allowed by the macroscopic data and
assumes that each occurs with equal  probability; this recipe
is suitable for an isolated system. An equivalent
characterization  can be given in terms of the canonical
ensemble, in which the
probability of a microscopic state is proportional to $\exp(-\beta E_s)$,
where $E_s$ is the
energy of
that state and $\beta$ is a parameter.
The principle of ``equivalence of ensembles" asserts that
these last two ensembles are
equivalent in the sense that, when the parameter $\beta$ is properly chosen
and the system has many degrees of freedom,
average  properties calculated in either ensemble are equal.

The parameter $\beta$ is generally called the ``inverse temperature" of the
system.
In many physical systems $T=1/\beta$ is indeed proportional to
 what one intuitively perceives to be the temperature, as it can
be gauged by
touching a system with one's finger. However, the parameter $\beta$ can be
viewed more abstractly,
as the parameter that makes the two ensembles be equivalent; in
incompressible turbulence, in which
the hydrodynamical motion and the vibrations of the molecules that make up
the fluid do not
affect each other, $\beta$ is unrelated to what
one normally calls the temperature of the fluid; $T$
does not have a readily intuitive interpretation;
standard thermodynamics
survives when one uses this less intuitive
notion of temperature (see e.g.~[5] Chap.~4). In a given
system, $\beta$ is a function of the energy kinetic energy $E$
of the turbulence and of whatever else is needed to describe
the system, just as in the kinetic theory of gases the
temperature is defined in terms of a kinetic energy of the
molecules. Our contention is that the intermediate  scales of
turbulence, those that are between the scales on which the
turbulence is stirred and those at which turbulence is
dissipated by viscosity, can be well-described by a theory
with the following properties: (i) These intermediate scales
are near a statistical equilibrium that is well-described as a
canonical or microcanonical equilibrium ensemble, (ii) this
equilibrium is characterized by values of macroscopical
parameters such as energy per unit mass or  squared  vorticity
per unit mass that define a critical state,
and (iii) the properties of this equilibrium can be derived from the
properties of
certain discrete systems by taking a special limit that we shall describe.

\section{Turbulence as a near-equilibrium process}

Turbulence {\it as a whole} cannot be described as  a
microcanonical or canonical equilibrium, for
if one stirs
a box full of fluid and then isolates the resulting flow, the outcome after a
long time is a state of rest;
an isolated system does not allow for the presence
of outside forces or an\linebreak
imposed shear
that would maintain the flow.
However, we apply canonical and microcanonical considerations only to
the range of
scales in turbulence intermediate between the dissipation scales and the
scale of the stirring, and the relevant question is whether the motion on
these
intermediate scales has enough time to settle down to an equilibrium on a time
scale in which one can
assume that little energy is added or subtracted from this intermediate
range,
so that one can view it as being approximately isolated;
in other words, we have to know
whether the
characteristic time of small-scale motion is short enough compared to the
characteristic time over which the flow as a whole changes appreciably.
According to the Kolmogorov
scaling in homogeneous isotropic turbulence, the structure
function $\<(u({\bold x}+{\bold r})-u({\bold x}))^2\>$,
where $u$ is the velocity along {\bf r},
is proportional to
$(\e r)^{2/3}$, where $r=|{\bold r}|$ and $\e$ is the rate of
energy dissipation per unit mass. Thus the characteristic time
(length/velocity) of an eddy of size $r$ is proportional
to $r/(\e r)^{1/3}=r^{2/3}\e^{-1/3}$.
The time scale over which the flow changes appreciably (in
the case of freely decaying turbulence, the characteristic time of over-all
decay) is $E/\e$, where
$E$ is the energy per unit mass. The ratio of these time scales is
proportional to $(\e
r)^{2/3}/E$, which tends to zero as $r\to 0$, and thus small
eddies (vortices) have time enough  time to approach a
statistical equilibrium.
In some of the turbulence literature, statistical equilibrium is
described by an equipartition ensemble, in which the energy
is equally distributed on the average between Fourier components
whose number tends to infinity. We have shown in earlier work
(see [6]) that this equipartition ensemble has no relevance to
fluid mechanics; its members have an infinite energy per unit
mass, and the corresponding  spectrum is proportional, in
three space dimensions, to $k^2$, where $k$ is a wave number.
In contrast, when we construct equilibria below,  we shall
take great care {\it not} to distribute energy evenly between
an increasing number of degrees of freedom. We shall instead
look for a sequence of finite-dimensional equilibria and a limit
designed so as to be meaningful
when the number of variables increases.

\section{Turbulence as a critical phenomenon}

Before discussing our construction in detail, a key observation must
be made regarding the relation between turbulence and critical phenomena.
We shall construct our ensemble and its probability density
as the limit of a sequence of ensembles of discrete vortex filaments living
on a lattice with mesh length $h$;
in each discrete ensemble
each  filament has a circulation $\Gamma$, The energy per unit
mass is $E$, and the total length of the filaments in the volume occupied
by a unit mass is $L$.
It is known [5] that, for a given choice of $\Gamma$ and $E$
and for each
mesh spacing $h$, the filaments break up into a
collection of
isolated vortex loops when  the dimensionless parameter
$\Lambda=L\sqrt{E}/\Gamma$ is small,
or  consist of a dense collection of smooth vortex lines
when
$\Lambda $ is large.
The transition between these cases, generated by varying either $E$, $L$,
or $\Gamma$  but keeping $h$ fixed,
has all the hallmarks of a phase transition in statistical mechanics: At
the points in the $\Gamma,E$
plane that correspond to the transition the fluctuations in various mean
quantities are large, the
derivative of $E$ with respect to the corresponding $\beta$ is large (i.e.,
a large change in $E$
changes the corresponding $\beta$ only a little), and the correlation
length (the length over which the
vortex configurations are correlated) is also large, as one expects from the
general theory of critical phenomena.
This phase
transition is related to the other, well-studied phase transitions in
statistical mechanics
[5].
We
claim that, to the extent that a lattice vortex system can be generated by a
fluid mechanical agency, the parameters $E,\Gamma,L$  must place it near this
phase transition.  We shall list now some obvious reasons, leaving some deeper
ones for the sequel.

(i) A key property of turbulence at large $Re$ is that it is subject to
scaling laws [1,2,3]. Scaling laws can hold for a many-particle or many-
vortex system
only near phase transition points; indeed, the existence of scaling
laws
is a general way of characterizing phase transition points. Phase
transition points where a scaling law holds are known as
``critical points".

(ii) Another salient feature of turbulence is that it dissipates energy. A
system
near equilibrium
dissipates  energy efficiently if the autocorrelation time of the fields that
describe it is long [8].
The
autocorrelation time
of a vorticity or velocity field is largest near ``critical" phase transition
points, and this is where energy dissipation can be substantial. A
version of this argument of particular relevance to vortex dynamics
is as follows:  A major mechanism of energy transfer from large to small
scales in three space dimensions is vortex intersection, which
reduces the scale of the vortical structures in the fluid. On the large
$\Lambda $ side
of the phase transition the vortex lines are smooth and there are few
vortex intersections, while on the small $\Lambda $ side of the phase
transition the vortices are small and isolated, and it is only on the
phase transition line that the vortices are both long enough and crumpled
enough to create a substantial number of intersections [5].

(iii) The hydrodynamical statistical equilibria have been discussed so far
without reference to how they may be produced by a fluid flow. In three
space dimensions vortex
lines stretch and fold; numerical calculations (see [5]) show that
vortex
stretching
pushes vortex systems to the neighborhood of the phase transition.
Furthermore,
a relatively simple calculation ([5], page 142) shows that, when the fluid
system is on the phase
transition line,  the energy spectrum
fluid system at the phase transition is the Kolmogorov spectrum.

\section{A discrete equilibrium model and its continuum limit}

We now describe in the discrete systems whose limit will be taken.
Consider vortex loops like the one in Figure 1, with
circulation $\Gamma$, the same for all the loops; these loops are elements
of impulse [4]; they  are located on a
$N\times N\times N$ mesh of mesh  size $h$, \ $Nh\!=\!1$.  An
incompressible velocity field is created by solving the
equations
\begin{equation}
 \Delta_h {\pmb\psi}=-{\pmb\xi} ,\qquad
{\bold u}=\nabla_h\times{\pmb\psi} \ ,
\end{equation}
where ${\pmb\xi}$ is the vorticity field, $\Delta_h$ is the
standard 6-point Laplacian in three space dimensions,
${\pmb\psi}$ is the vector potential, $\nabla_h$ is the
differencing operator implemented with forward differences;
equation (1) is the discrete analog of the
usual definition of the velocity in terms of the vorticity; the domain is
assumed
to be periodic with period $1$. The kinetic energy of the system is defined as
\begin{equation}
E_h=\textstyle{\frac 12}\Sigma(u_1^2+u_2^2+u_3^2)h^3 \ ,
\end{equation}
where $u_1,u_2,u_3$ are the components of ${\bold u}_h$
and the summation is over all the nodes in the lattice.
The energy $E_h$ being fixed, the appropriate equilibrium ensemble is the
microcanonical ensemble.
For a given $h$, the various configurations of the equilibrium ensemble
can be
sampled by ``microcanonical sampling" [7].
The parameter $\beta$ in the equivalent
canonical ensemble can be determined in the course of calculating averages.

If one decreases $h$ while keeping the energy $E_h$
fixed
the ``temperature" $T=1/\beta$ of the system decreases; heuristically, $T$ is
the energy per degree of freedom, and as the number of degrees of
freedom increases, there is less energy for any one of them -- this is exactly
the situation in the ``equipartition" ensemble discussed above.
The limiting ensemble does not exist. On the other hand,
suppose one decreases $h$ while keeping fixed not only the
energy but also some measure of the total vorticity -- for
example, the sum
$Z_2=h^3\sum|{\pmb\xi}|^2$, where ${\pmb\xi}$ is the vorticity field formed
by the sides of the vortex loops described above and in Figure 1,
and the summation is over all the nodes in the lattice.
$Z_2$ cannot be viewed as an accurate approximation of
the enstrophy (we  would have to know how the vorticity is
distributed   across each leg of
the loops),
but it is a sensitive gauge of the total vorticity in the
system.
If one decreases $h$ while keeping both $E_h$ and $Z_2$ constant, as one
can easily do with microcanonical sampling, the temperature of the system
tends to infinity and beyond; indeed, temperatures beyond infinity, the
so-called
``negative temperatures", are typical of two dimensional turbulence [5]
where additional integrals of motion constrain the vorticity. The
heuristic explanation of
this phenomenon is that in the presence of two constraints there  are ever
fewer states that satisfy both constraints; the effective number of degrees
of freedom
tends to one and the temperature increases. In the limit $h\to0$ one finds a
probability measure concentrated on a single organized structure, which is
reasonable in two dimensions but not in three. Thus a fixed bound on $Z_2$
produces
too much organization while an absence of bounds on the vorticity produces too
little organization; the ensemble we seek should be between those two extremes.

What we shall do is pick some finite value of $\beta$ and
adjust $Z_2$ so that the system
assume this inverse temperature $\beta$ for each $h$. We already showed
elsewhere [6]
that this can be done in two
space dimensions. The automatic way of finding the appropriate value
of $Z_2$ described in [6]
turned out to be impractical in three space dimensions, and we proceed by
simple tabulation:
Pick various bounds for $Z_2$, calculate the corresponding values of
$\beta$ and home in on the value
of $\beta$ that one wants to have by successive improvement guided by a human.
The ensemble produced in this
way is intermittent -- the vorticity occupies only a fraction
of the available sites, and has
hidden coherent structures- the temperature being finite, there is
effectively a finite number
of degrees of freedom, fixed as $h\to0$, but their nature and location is
not explicitly known.
The finiteness of $\beta$ ensures that one can define a probability density
in the limit state.

The choice of the fixed $\beta$ is arbitrary, but as $h\to0$ the physical
results
become independent of that fixed $\beta$ because the
resulting systems converge
to the critical line: A refinement of the mesh is tantamount to a
change in the scale of the vortex configurations; the critical points
are the loci of the states that are invariant under a change in scale,
and attract all other states as $h\to0$, as long as their temperature is
finite and non-zero.
Controlling the rate of growth of the enstrophy as $h\to0$ is
consistent
with the physics of three-dimensional turbulence:
Turbulence is intermittent, and vorticity concentrates on small scales.
As $h$ is decreased ever more vorticity is revealed, at a rate that determines
the vorticity spectrum  and hence the  energy spectrum.

\section{Some numerical results}

We now display some results obtained in
calculations carried out as was just explained.
The lattices we could afford so far
are small, and the results are not converged; indeed,
if one views the mesh size $h$ as an analog of a viscous length, the
scaling
derived in [1,2,3] shows that the effect of the mesh size
decreases like
$\frac{1}{\ln |h|}$
; this is the rate of convergence of analogous computations
in two space
dimensions as well [6], and even in two dimensions the calculations
cannot be brought to convergence without relying on
explicit renormalization
group parameter flow which we do not have in three dimensions. We
expect here
only qualitative results.
Analogous, more complete computations in simpler models have been displayed
elsewhere.

In Figure 2 we display
the variation of the computed $\beta=T^{-1}$ as a function of the lattice
size $N$
for a fixed energy $E_h=100$,
with no restrictions on the ``enstrophy" $Z_2$. The graph shows
that $\beta$ is proportional
to $N^3$, as one can expect from equipartition; without constraints on the
vorticity and at a fixed, finite energy,
the temperature of the system tends to zero.
In Figure 3 we exhibit the values of $\beta$ obtained by bounding
$Z_2$ with $N=16$: We first found the value $Z_{2,0}$ of $Z_2$
that corresponded to equilibrium with no bound on $Z_2$, and
then repeated the calculation requiring that $Z_2\le
Z_{20}-\delta Z_2$ for various values of
$\delta Z_2$. This graph shows that bounding $Z_2$ does indeed
keep the ``temperature" $T=\beta^{-1}$ away from zero.

In Figure 3  we display the shape of the second order
structure function $ R(r)=C\<(u({\bold x}+{\bold r})-u({\bold x}))^2\>$ as a
function of the separation $r=|{\bold r}|$, with $N=16$,
with $u$ the velocity  component in the direction of {\bf r}, at two values
of
$\beta$:
$\beta=40.2$ and $\beta=3$ obtained with differing values of $\delta
Z_2$.
The constant $C$ equals $(4/3)/ \<(u({\bold x}))^2\>$ so as to
make
$R(r)=1$ for large $r$.
While no convergence can be expected on this crude grid,
the graph does show that controlling $\beta$ produces more plausible
structure functions, in particular structure functions that do not tend to
a $\delta$ function, as one would get in an equipartition ensemble.
The number of Monte-Carlo steps needed for convergence is large, and a finer
grid could not be used in the present state of our algorithm.

One point we can hope to clarify with the model we have is whether the
higher moments
of the velocity field remain bounded -- a major issue in the
scaling of the local structure
of turbulence. We have deduced in [1,2,3], from scaling and heuristic
arguments, that the moments of order
$p\le 3$
of the velocity field (and in the case of wall-bounded turbulence,
the mean of the velocity gradient) can be obtained by an expansion
about a vanishing viscosity limit (in powers of $\frac{1}{\log Re}$,
where $Re$ is the Reynolds number. We conjectured that the reason
this apparently could not be done for higher moments was that the
higher order moments, of order $p\ge 4$, diverged as the viscosity
tended to zero. Heuristic argument why this should be so were given
in [1,2,3,6]. If this conjecture is true, the Kolmogorov scaling is
the local structure of turbulence holds as long as the vanishing
viscosity limit exists. The calculations described here, with $N=16$
and $E_h=100$, give for $\<u_1^2\>$, where $u_1$ is a component
of the velocity {\bf u}, the value $\<u_1^2\>=66$, (a value
dictated by the imposed $E_h$); however,
$\<u_1^4\>=1.3\cdot10^5$, $\<u_1^6\>=4.2\cdot10^7$, with a
standard deviation of 10\%. The odd moments are zero up to a few standard
deviations;
however, one should note that for the
fifth moment the standard deviation is large, of order
$10^5$, reflecting the very large value of the tenth
moment. These values change little with lattice size $N$,
as one expects from the analysis in [1,2,3] which shows
that a significant change in the moments requires a
significant change in $\log N$. (Note that an analogous argument
applies to experimental data; a variation of the higher moments as
the viscosity changes may be very hard to detect).
There preliminary results are consistent with
our earlier conjectures.

\section{Acknowledgements}
We would like to thank Prof.~P.~Colella,  Dr.~A.~Kast and
Dr.~R.~Kupferman for helpful discussions and comments.

\newpage
\begin{center}
{\bf\large References}
\end{center}

\begin{enumerate}

\item
G.I. Barenblatt and A.J. Chorin, Scaling laws and vanishing viscosity
limits for wall-bounded shear flows and for local structure in developed
turbulence, {\it Comm. Pure Appl. Math.} {\bf 50}, 381--391 (1997).

\item
G.I. Barenblatt and A.J. Chorin, New perspectives in turbulence: Scaling laws,
asymptotics and intermittency, {\it SIAM Review} {\bf 40}, 265--291 (1998)

\item
G.I. Barenblatt, A.J. Chorin, and V.M. Prostokishin, Scaling laws in
fully developed
turbulent pipe flow, {\it Appl. Mech. Rev.}, {\bf 50}, 413--429, (1997).

\item
T. Buttke and A. Chorin, Turbulence calculations in magnetization
variables, {\it Appl. Num. Math.} {\bf 12}, 47--54 (1993)

\item
A. J. Chorin, {\it Vorticity and Turbulence}, Springer, 1994.

\item
A.J. Chorin, New perspectives in turbulence, {\it Quart. Appl. Math.,} 1998,
in press.

\item
H. Gould and J. Tobochnik, {\it An Introduction to Computer Simulation Methods,
Applications to Physical Systems}, part 2, Addison-Wesley, Reading, 1988.

\item
L. Landau and E. Lifshitz, {\it Statistical Physics}, Vol. 1, Pergamon Press,
Oxford, 1980,
pp. 368--389.

\end{enumerate}

\newpage
\centerline{\bf Figure Captions}

\noindent Figure 1. A basic vortex loop.

\bigskip\noindent
Figure 2. The variation of the inverse temperature $\beta$ with $N$
in the absence of vorticity constraints.

\bigskip\noindent
Figure 3.
The increase of the temperature as the bound
on the enstrophy is decreased (The notation is explained in the text).

\bigskip\noindent
Figure 4. Scaled structure function
$R(r)=C\<(u({\bold x}+{\bold r})-u({\bold x}))^2\>$,
$C=4/(3
\<u^2({\bold x})\>$, $r=|{\bold r}|$,
for  $N=16$, with $\beta=3 $ and with $\beta=40.3$
(Small $\beta$ corresponds to a
large temperature $T$ and requires a small
bound on the enstrophy $Z_2$).

\end{document}